\documentclass[10pt]{amsart}

\usepackage{graphicx}

\usepackage[USenglish]{babel}
\usepackage{latexsym}
\usepackage{amsmath}
\usepackage{amsfonts}
\usepackage{amssymb}
\usepackage{esint}
\usepackage[all]{xy}

\newcommand{\D }{\Delta }

\newcommand{\Sg }{\Sigma}

\newcommand{\be}{\begin{equation}}
\newcommand{\ee}{\end{equation}}
\newenvironment{pf}{\noindent{\sc Proof}.\enspace}{\rule{2mm}{2mm}\medskip}

\newcommand{\R}{\mathbb{R}}

\newcommand{\N}{\mathbb{N}}
\newcommand{\no}{\noindent}

\newtheorem{theorem}{Theorem}[section]
\newtheorem{proposition}[theorem]{Proposition}
\newtheorem{example}[theorem]{Example}

\newcommand{\bpr}{\begin{proposition}}
\newcommand{\epr}{\end{proposition}}
\newcommand{\bex}{\begin{example}\rm}
\newcommand{\eex}{\end{example}}

\begin{document}

\newtheorem{lem}{Lemma}[section]
\newtheorem{pro}[lem]{Proposition}
\newtheorem{thm}[lem]{Theorem}
\newtheorem{rem}[lem]{Remark}
\newtheorem{cor}[lem]{Corollary}
\newtheorem{df}[lem]{Definition}

\title[The Mean Field equation on compact surfaces]{New existence results for the mean field equation on compact surfaces via degree theory}

\author {Aleks Jevnikar}

\address{SISSA, via Bonomea 265, 34136 Trieste (Italy).}

\email{ajevnika@sissa.it}

\keywords{Geometric PDEs, Leray-Schauder degree, Mean field equation.}

\subjclass[2000]{ 35J20, 35J61, 35R01.}

\begin{abstract}
We consider the following class of equations with exponential nonlinearities on a closed surface $\Sg$:
$$
  - \D u = \rho_1 \left( \frac{h \,e^{u}}{\int_\Sg
      h \,e^{u} \,dV_g} - \frac{1}{|\Sg|} \right) - \rho_2 \left( \frac{h \,e^{-u}}{\int_\Sg
      h \,e^{-u} \,dV_g} - \frac{1}{|\Sg|} \right),
$$
which arises as the mean field equation of the equilibrium turbulence with arbitrarily signed vortices. Here $h$ is a smooth positive
function and $\rho_1, \rho_2$ two positive parameters. 
By considering the parity of the Leray-Schauder degree associated
to the problem, we prove solvability for $\rho_i \in (8\pi k, 8\pi(k+1)),\, k \in \N$. Our theorem provides a new
existence result in the case when the underlying manifold is a sphere and gives a completely new
proof for other known results.
\end{abstract}

\maketitle

\section{Introduction}

\no We are concerned with the equation
\begin{equation}
-\D u = \rho_1 \left(\frac{h\,e^u}{\int_\Sigma h\,e^u \,dV_g} - \frac{1}{|\Sigma|}\right) - \rho_2 \left(\frac{h\,e^{-u}}{\int_\Sigma h\,e^{-u} \,dV_g} - \frac{1}{|\Sigma|} \right) \hspace{0.3cm} \mbox{on $\Sigma$}, \label{eq}
\end{equation}
where $\D = \D_g$ is the Laplace-Beltrami operator, $\rho_1, \rho_2$ are two non-negative parameters, $h:\Sigma \rightarrow \mathbb{R}$ is a smooth positive function and $\Sigma$ is a compact orientable surface without boundary with Riemannian metric $g$ and total volume $|\Sigma|$. For the sake of simplicity, we will assume throughout this work that
$|\Sigma|=1$, which is no loss of generality thanks to a trivial rescaling argument'.

Problem \eqref{eq} plays an important role in mathematical physics as a mean field equation for the equilibrium turbulence with
arbitrarily signed vortices. It was first obtained by Joyce and Montgomery \cite{joy-mont} and  by Pointin and Lundgren \cite{point-lund} by means of different statistical arguments. Later, several authors adopted this model; we refer for example to \cite{cho}, \cite{lio}, \cite{mar-pul}, \cite{new} and the references therein. The case $\rho_1 = \rho_2$ is also related to the study of constant mean curvature surfaces, see \cite{wen1}, \cite{wen2}.

Before describing the main features of the problem and the known results, let us first consider the case $\rho_2 = 0$, namely the following standard Liouville-type equation:
\begin{equation} \label{eq2}
  - \D u = \rho \left( \frac{h \,e^{u}}{\int_\Sg
      h \,e^{u} \,dV_g} - 1 \right).
\end{equation}
Equation \eqref{eq2} appears in conformal geometry in the problem of finding a conformal metric for which the Gauss curvature is a prescribed function on $\Sg$, see \cite{bah-cor}, \cite{cha}, \cite{cha2}, \cite{li}, \cite{scho-zha}. Indeed, setting $\tilde g= e^{2v}g$, the Laplace-Beltrami operator for
the deformed metric is given by $\D_{\tilde g} = e^{-2v}\D_g$ and the change of the Gauss curvature is ruled by
$$
	-\D_g v = K_{\tilde g} e^{2v} - K_g,
$$  
where $K_g$ and $K_{\tilde g}$ are the Gauss curvatures of $(\Sg, g)$ and of $(\Sg, \tilde g)$ respectively.
Problem \eqref{eq2} also arises in mathematical physics as a mean field equation of Euler flows, see \cite{cal}, \cite{kies}. The literature on \eqref{eq2} is broad, and there are many results regarding existence, blow-up analysis, compactness of solutions, etc, see \cite{djlw}, \cite{djadli}, \cite{mal}, \cite{tar}.

\medskip

As many geometric problems, also \eqref{eq2} carries a lack of compactness, as its solutions might blow-up. It was proved in \cite{bre-mer}, \cite{li1} and \cite{li-shaf} that a quantization phenomenon occurs in this case. More precisely, taking a blow-up point $p$ for a sequence $(u_n)_n$ of solutions, we have
\begin{equation} \label{quant}
\lim_{r \to 0} \lim_{n \to + \infty}  \int_{B_r(p)} h \, e^{u_n} dV_g = 8 \pi.
\end{equation}
Roughly speaking, each blow-up point carries a quantized local mass. Moreover, the limit profile of solutions becomes close to a {\em bubble}, namely a function $U_{\lambda,p}$ defined as 
$$
	U_{\lambda,p}(y) = \log \left( \frac{4\lambda}{\bigr( 1 + \lambda \,d(p,y)^2 \bigr)^2} \right),
$$
where $y\!\in\!\Sg,\, d(p,y)$ stands for the geodesic distance and $\lambda$ is a large parameter. This limit function can be viewed as the logarithm of the conformal factor of the stereographic projection from $S^2$ onto $\R^2$, composed with a dilation.

Combining the local quantization \eqref{quant} with some further analysis, see for example \cite{bat-man}, \cite{bre-mer}, we have that the set of solutions
to \eqref{eq2} is uniformly bounded in $C^{2,\alpha}$, for any fixed $\alpha \in (0,1)$, provided $\rho \notin 8\pi\N$. It follows that one can define the Leray-Schauder degree associated to problem \eqref{eq2} with $\rho\in(8k\pi, 8(k+1)\pi), k\in\N$. In \cite{li1} it was shown that the degree is $1$ when $\rho < 8\pi$. By the homotopic invariance of the degree, it is easy to see that the same is independent of the function $h$, the metric of $\Sg$ and it is constant on each interval $(8k\pi, 8(k+1)\pi)$. In fact it depends only on $k\in\N$ and the topological structure of $\Sg$, as was proved in \cite{chen-lin}, where the authors provide the degree-counting formula
\begin{equation} \label{deg}
	\mbox{deg}(\rho) = \frac{1}{k!}(-\chi(\Sg) + 1) \cdots (-\chi(\Sg) + k),
\end{equation}
where $\chi(\Sg)$ denotes the Euler characteristic of $\Sg$. The proof of this result is carried out by analyzing the jump values of the degree after $\rho$ crossing the critical thresholds. Later, this result was rephrased in \cite{mal1} with a Morse theory point of view.

\medskip

On the other hand, in the general case when $\rho_2 \neq 0$, namely for equation \eqref{eq}, there are fewer results: for example the refined blow-up analysis of solutions of \eqref{eq} is not yet fully developed. Nevertheless, it was proven in \cite{jwyz} that the blow-up phenomenon yields a quantization property; for a blow-up point $p$ and a sequence $(u_n)_n$ of solutions relatively to $(\rho_{1,n}, \rho_{2,n})$ the authors obtained
$$
    \lim_{r \to 0} \lim_{n \to + \infty} \rho_{1,n} \frac{ \int_{B_r(p)} h \, e^{u_n} \,dV_g }{ \int_\Sg h \, e^{u_n} \,dV_g } \in 8 \pi \N, \qquad
    \lim_{r \to 0} \lim_{n \to + \infty} \rho_{2,n} \frac{ \int_{B_r(p)} h \, e^{-u_n} \,dV_g }{ \int_\Sg h \, e^{-u_n} \,dV_g } \in 8 \pi \N.
$$
Moreover, the case of positive multiples of $8 \pi$ may indeed occur, see \cite{es-we} and \cite{gr-pi}. The latter volume quantization implies that the set of solutions is compact for $(\rho_1, \rho_2) \notin (8 \pi \N \times \R) \cup (\R \cup 8 \pi \N)$. 

It follows that, as before, the degree associated to \eqref{eq} can still be defined outside this set of parameters. However, this strategy has not been yet investigated and the existence results mostly rely on a variational approach.

\medskip

Let us now briefly discuss the history of the problem and the previous results. First of all, the case $\rho_1, \rho_2 < 8\pi$ was considered in \cite{oh-su}. Ohtsuka and Suzuki proved that in this regime the associated energy functional is bounded
from below and coercive and therefore solutions can be found as global minima. If $\rho_i > 8\pi$ for some $i=1,2,$ then the problem becomes subtler as the functional is not bounded and a minimization technique is not possible any more. The first result in this direction is given in \cite{jwyz} for $\rho_1 \in (8\pi, 16\pi)$ and $\rho_2<8\pi$. Via a blow-up analysis Jost et al. proved the existence of solutions on a smooth, bounded, non-simply connected domain in $\R^2$ with homogeneous Dirichlet boundary condition. Later, in \cite{zhou}, Zhou generalized this result to any compact surface without boundary by reducing the problem to the one with $\rho_2 = 0$ and applying min-max theory. 

The next step was made in \cite{jev}, where for the first time a doubly supercritical regime is considered, namely $\rho_i \in (8\pi, 16\pi)$ for $i=1,2.$ The existence of solutions is proved using a min-max scheme based on a detailed description of the low sublevels of the energy functional. Finally, this was generalized in \cite{todatori} for a generic choice of the parameters $\rho_1 \in (8k\pi, 8(k+1)\pi)$ and $\rho_2 \in (8l\pi, 8(l+1)\pi),\, k,l \in \N$, under the assumption that $\Sg$ has positive genus.

\section{The main result}

We attack the problem with a different point of view and for the first time we analyze the associated Leray-Schauder degree. This is done in the spirit of \cite{mrd}, where a Toda system of Liouville equations arising from Chern-Simons theory was analyzed. More precisely, we study its parity and we observe that when both parameters stay in the same interval, i.e. $\rho_i \in (8k\pi, 8(k+1)\pi),\, k \in \N$ for $i=1,2$, the degree is always odd. The main result is the following.

\begin{thm}\label{teo}
Let $h>0$ be a smooth function and suppose \mbox{$\rho_i \in (8k\pi, 8(k+1)\pi)$,} $k \in \N$ for $i=1,2$. Then problem \eqref{eq} has a solution.
\end{thm}

Observe that we recover the result of \cite{jev} and some cases of \cite{todatori}: when $\Sg$ is homeomorphic to $S^2$ the above theorem yields a new existence result.

\

\begin{pf}
For some $\alpha \in (0,1)$ let $C^{2,\alpha}_0 (\Sg)$ be the class of $C^{2,\alpha}$ functions with zero average. Consider now the mapping $T : C^{2,\alpha}_0 (\Sg) \to C^{2,\alpha}_0 (\Sg)$ defined by
\begin{equation} \label{map}
	T(u) = (-\D)^{-1} \left( \rho_1 \left( \frac{h \,e^{u}}{\int_\Sg
      h \,e^{u} \,dV_g} - 1 \right) - \rho_2 \left( \frac{h \,e^{-u}}{\int_\Sg
      h \,e^{-u} \,dV_g} - 1 \right) \right),
\end{equation}
where $(-\D)^{-1}f, f\in C^\alpha(\Sg)$, is intended as the solution $v$, with zero average, of the problem $-\D v = f$, which is unique. We are concerned with the map $\Psi = Id - T$ and the solutions of equation \eqref{eq} will correspond to zeros of $\Psi$. 

Clearly, by elliptic regularity theory the operator $T$ is compact. Moreover, the set of the solutions is compact for parameters $(\rho_1, \rho_2) \notin (8 \pi \N \times \R) \cup (\R \times 8 \pi \N)$. Therefore, we can consider the associated degree deg$\left(\Psi_{(\rho_1,\rho_2)},B_r(0),0\right)$ which is well-defined for $r$ sufficiently large.

Consider now $\rho_i \in (8k\pi, 8(k+1)\pi),\, k \in \N$ for $i=1,2$. Letting $\rho = \frac 12 (\rho_1 + \rho_2)$, we perform the following homotopy which
takes place in a connected component of $\R^2 \setminus ((8 \pi \N \times \R) \cup (\R \times 8 \pi \N))$:
$$
	\mathfrak h(t) = (1-t)(\rho_1,\rho_2) + t (\rho,\rho).
$$
From the fact that the degree is constant along homotopies we obtain that
$$
	\mbox{deg}\left(\Psi_{(\rho_1,\rho_2)},B_r(0),0\right) = \mbox{deg}\left(\Psi_{(\rho,\rho)},B_r(0),0\right).
$$
Observe now that by the structure of $T$ we deduce
$$
	\Psi_{(\rho,\rho)}(-u) = -\Psi_{(\rho,\rho)}(u).
$$
Therefore, we conclude that $\Psi_{(\rho,\rho)}$ is an odd operator. By the Borsuk theorem, see \cite{kyner}, it follows that the associated degree is odd an hence non zero. This guarantees us the existence of a solution to equation \eqref{eq}.
\end{pf}

\

\begin{center}
\textbf{Acknowledgements}
\end{center}

\no The author would like to thank Professor Andrea Malchiodi for his kind help in preparing this paper.

The author is supported by the PRIN project \emph{Variational and perturbative aspects of nonlinear differential problems}.

\end{document}